\documentclass[12pt]{amsart}
\usepackage[dvipsnames]{xcolor}
\usepackage[top= 1in, bottom=  1in,left=1in,right=1in]{geometry}

\usepackage{latexsym, amsmath, mathtools, amssymb,amsthm,amsopn,amsfonts, eucal, dsfont, bbm,float}

\definecolor{Awesome}{rgb}{1.0, 0.13,0.32}
\definecolor{BleuDeFrance}{rgb}{0.19, 0.55,0.91}
\newcommand{\R}[1]{\mathbb{R}^{#1}}

\newcommand{\C}{\mathbb{C}}

\newcommand{\A}{\mathcal{A}}
\newcommand{\M}{\mathcal{M}}
\renewcommand{\L}{\mathbf{L}}
\newcommand{\e}{\varepsilon}
\newcommand{\LL}{\mathcal{L}}

\newcommand{\loc}{\mathrm{loc}}
\def\mR{\mathbb{R}}
\def\mC{\mathbb{C}}
\def\mZ{\mathbb{Z}}
\newcommand{\eps}{\varepsilon}
\newcommand{\abs}[1]{\left\vert{#1}\right\vert}
\DeclareMathOperator{\dist}{dist}
\DeclareMathOperator{\vol}{vol}
\newtheorem{thm}{Theorem}
\newtheorem{prop}{Proposition}
\newtheorem{definition}{Definition}
\newtheorem{rmk}{Remark}
\newtheorem{lem}{Lemma}
\newtheorem{cor}{Corollary}

\author{Dmitry Golovaty}
\address{Department of Mathematics, The University of Akron, Akron OH 44325}
\email{dmitry@uakron.edu}

\author{Alberto Montero}
\address{Departamento de Matem\'atica y Ciencia de la Computaci\'on, Universidad de Santiago de Chile}
\email{alberto.montero@usach.cl}

\author{Etienne Sandier}
\address{D\`epartement de Math\'ematiques
Universit\'e Paris 12 
94010  Cr\'eteil Cedex, France
}
\email{sandier@u-pec.fr}

\author{Peter Sternberg}
\address{Department of Mathematics, Indiana University,
Bloomington, IN 47405}
\email{sternber@iu.edu}

\title{On sections of complex line bundles over surfaces minimizing a Ginzburg-Landau energy}

\begin{document}

\begin{abstract}
    In this work we extend some of the results of \cite{IgnatJerrard} for Ginzburg-Landau vortices of tangent vector fields on 2-dimensional Riemannian manifolds to the setting of complex hermitian line bundles. In particular, we elucidate the locations of vortices for the cases of $Q$-tensors and their higher-rank analogs on a sphere. 
\end{abstract}

\maketitle

\section{Introduction}
In this note we extend the asymptotic analysis of the minimization of the Ginzburg-Landau energy to the setting of a complex hermitian line bundle $\L$ over a closed two-dimensional surface $\M$. 

When both the connection $1$-form and the section of the bundle are considered as unknowns, this is known as the abelian Higgs model (see \cite{JaffeTaubes}, and \cite{BethuelRiviere}, \cite{Orlandi} for the asymptotic analysis of minimizers).

The minimization with respect to the section only, the connection $1$-form being fixed,  was considered more recently. In \cite{IgnatJerrard}, Ignat and Jerrard  analyze this for the case where $\L$ is the tangent bundle with its Levi-Civita connection, so that a section of $\L$ is a tangent vector field on $\M$. Their analysis carries over straightforwardly to a general bundle equipped with a connection,  as we show. 

Our motivation lies in the context of liquid crystal modeling where one often wishes to capture nonorientable configurations, leading to the consideration of a Ginzburg-Landau type energy defined not on vector fields but rather on so-called $Q$-tensor bundles. In fact, here we consider a more general theory that encompasses higher order symmetric traceless tensor bundles, as we review in Section \ref{sec:high}. These objects have applications in meshing and physics and are discussed in more detail in \cite{MR4572144,MR4315483,10.1145/3366786}.

Specifically, the energy to be studied here is given by
\begin{equation}
E_\e[u]=\frac12\int_{\M}|\nabla_Au|^2+\frac{1}{2\varepsilon^2}{(1-|u|^2)}^2,\label{GLenergy}
\end{equation}
where $\e>0$, $A$ is a connection on $\L$ and $u$ is a section. In Section \ref{linebundles} we provide background on line bundles, connections and sections, but here we point out that locally, such a section $u$ is simply given by its trivialization, which takes the form of a map from an open patch on $\M$ to the complex plane $\C$. 

Going back to the seminal work of Brezis, Bethuel and Helein \cite{BBH}, the central goal in studying minimizers of $E_\e$ is to determine the number and asymptotic location of the singularities that for topological reasons are forced to exist. In the present setting, singularities necessarily emerge when the Euler number of $\L$ is non-zero. In the context of liquid crystals these singularities represent defects whose appearance represent key features of the liquid crystal texture on the surface. As in \cite{BBH} and \cite{IgnatJerrard}, we will characterize the asymptotic location of these singularities as those points on $\M$ minimizing a certain renormalized energy $W$ arising as the $O(1)$ term, following a logarithmic leading term, in an expansion for the minimal energy. An explicit formula for $W$ is given in Proposition \ref{explicit}.

Our main result is the following:
\begin{thm}
\label{thm:1}
Assume that $\L$ is a complex hermitian line bundle over a closed surface $\M$ and $A$ is a connection on $\L$. 
Let $u_\varepsilon$ be a minimizer of $E_\e.$ Then there exist points $a_1,\ldots,a_d$ where $d$ is the absolute value of the Euler number $e(\L)$ such that $u_\e$ converges to $u_*$ in $C^k_{\mathrm{loc}}(\M\backslash\{a_1,\ldots,a_d\})$
as $\eps\to0$ for any $k$, where $u_*$ is a {\em canonical harmonic section}  with singularities of degree one about the points $a_1,\ldots,a_d$ (see definition below).  Furthermore, $E_\eps[u_\eps]=\pi d\log{\frac{1}{\eps}}+d\gamma+W(a_1,\ldots,a_d)+o(1).$ Here, for any collection of points $b_1,\dots,b_d$ on $\M$ we let  $W(b_1,\ldots,b_d)$ be the minimum  of \begin{equation}
    \tilde{E}[v,b_1,\ldots,b_d]:=\lim_{\rho\to0}\left[\frac{1}{2}\int_{\M\backslash\cup_{j=1}^d B_\rho(b_j)}|\nabla_A v|^2+\pi d\log{\rho}\right]\label{Wdefn}
    \end{equation}
over all canonical harmonic sections $v$  with singularities of degree one about the points $b_1,\ldots,b_d$.

The points $(a_1,\ldots,a_d)$ minimize $W$ over $\M^d$ and the section $u_*$ minimizes $\tilde{E}[\cdot,a_1,\ldots,a_d]$ over canonical harmonic sections singular at $a_1,\ldots,a_d$.
\end{thm}

\begin{rmk}
In \cite{IgnatJerrard} there is 
a proof of $\Gamma$-convergence but here we focus only on the limiting behavior of minimizers. This simplifies the proof of the result greatly. However, we do not see any impediment to establishing such a more general result in our setting.
\end{rmk}

A more explicit expression for $u_*$ leading to an explicit formula for $W$ in terms of the Green's function for $\M$ was derived in \cite{IgnatJerrard} for the specific case of the tangent bundle on $\M$. We adapt this procedure to the present, more general case in Section \ref{sec:W}. We then use the general formula for $W$ to specialize the results to bundles of higher rank symmetric traceless tensors in the next theorem. The relevant background is provided in Section \ref{sec:high}.
\begin{thm}\label{tetra}
    If $\M=\mathbb{S}^2$ and $\L$ is the bundle $T^k_{sym, tr}\mathbb S^2$ of symmetric traceless rank $k$ tensors on $\mathbb{S}^2,$ then the number of singularities of the limiting canonical harmonic section is $d=2k$. Moreover, the location of these singularities can be characterized as follows: 
    \begin{itemize}
        \item If $k=2$ then $(a_1,\ldots,a_4)$ lie at the vertices of a regular tetrahedron.
        \item If $k=3$ then $(a_1,\ldots,a_6)$ lie at the vertices of the cross-polytope.
        \item If $k=6$ then $(a_1,\ldots,a_{12})$ lie at the vertices of a regular icosahedron.
    \end{itemize}
\end{thm}

\begin{rmk}
    When $k=2,$ the bundle is precisely that of $Q$-tensors \cite{MR2893657}; thus singularities of a uniaxial nematic liquid crystal on the surface of a sphere form a perfect tetrahedron, as is well known from experimental observations.
\end{rmk}

\section{Background on line bundles}\label{linebundles}
In this section, for the convenience of the reader, we provide some basic background on hermitian line bundles over surfaces, cf. \cite{BottTu, Hall}. To this end, let $\M$ be a surface. We let $\{U_\alpha\}_\alpha$ be a {\em simple} cover of $\M,$ i.e., a cover of $\M$ by simply connected open patches such that any intersection $U_{\alpha\beta} = U_\alpha\cap U_\beta$ or  $U_{\alpha\beta\gamma} = U_\alpha\cap U_\beta\cap U_\gamma$ is simply connected as well. Then  a complex line bundle $\L$ over $\M$ is defined by transition functions $\varphi_{\alpha\beta}:U_{\alpha\beta}\to\mR $ satisfying the condition that $\varphi_{\alpha\beta} + \varphi_{\beta\gamma} + \varphi_{\gamma\alpha}$ defined on $U_{\alpha\beta\gamma}$ is $2\pi\mZ$-valued.

A {\em section} $u$ of $\L$ is defined by its trivializations 
$\{u_\alpha\}_\alpha$
which are maps $$u_\alpha:U_\alpha\to \mC $$ such that $u_\alpha = e^{i\varphi_{\alpha\beta}} u_\beta$ on $U_{\alpha\beta}$. If we further choose local coordinates on $U_\alpha$, then we can view $u_\alpha$ as a map from an open subset of $\mR^2$ to $\mC$.
Notice that we can define an inner product between two sections $u$ and $v$ using the complex structure of the trivializations on $U_\alpha$, i.e.
$\langle u_\alpha,v_\alpha\rangle=\bar{v}_\alpha u_\alpha$, and this inner product is independent of $\alpha$.

A {\em connection} $A$ on $\L$ is defined by a set of real-valued 1-forms $\{A_\alpha\}_\alpha$, where $A_\alpha$ is defined on $U_\alpha$ and where 
\begin{equation}
    \label{eq:conn}
    A_\alpha = A_\beta + d\varphi_{\alpha\beta}
\end{equation}
on $U_{\alpha\beta}$.

The {\em covariant derivative} of a section $u$ relative to the connection $A$ in the direction of the tangent field $\tau$ is a section denoted $\partial_\tau^A u$ and defined by its trivializations
\begin{equation}\label{covdev}(\partial_\tau^A u)_\alpha = \tau\cdot u_\alpha - i A_\alpha(\tau) u_\alpha,\end{equation}
where $\tau\cdot u_\alpha$ is the derivative of the complex-valued function $u_\alpha$ in the direction $\tau$. Note that we have
\begin{equation*}
    \begin{split}
  (\partial_\tau^A u)_\alpha &= \tau\cdot (e^{i\varphi_{\alpha\beta}}u_\beta) -  i \left(A_\beta + d\varphi_{\alpha\beta}\right)(\tau) e^{i\varphi_{\alpha\beta}}u_\beta\\
  &= e^{i\varphi_{\alpha\beta}}\left( \tau\cdot u_\beta - iA_\beta(\tau) u_\beta + id\varphi_{\alpha\beta}(\tau) u_\beta - id\varphi_{\alpha\beta}(\tau) u_\beta\right)\\
  &= e^{i\varphi_{\alpha\beta}} (\partial_\tau^A u)_\beta,\\
\end{split}
\end{equation*}
so that this indeed defines a section of the bundle.
We denote by $d_A u$ the map $\tau \to \partial_\tau^A u$.

The two-form $dA$ is globally defined since $dA_\alpha = dA_\beta$ on $U_{\alpha\beta}$. It is called the {\em curvature} 2-form and $\kappa = *dA$ is the {\em curvature} of the connection $A$.

Next, given points $\{b_1,\ldots,b_d\}\subset\M,$ we specialize to the case of a smooth section $v$ of modulus $1$ defined on $\M\setminus\{b_1,\ldots,b_d\}$. It has trivializations $v_\alpha,$ where $v_\alpha:U_\alpha\backslash\{b_1,\ldots,b_d\}\to\mathbb{S}^1.$ If we denote by $\theta_\alpha$ a local lifting of $v_\alpha$ then the covariant derivative of $v$, as given by \eqref{covdev} can be written in the form
\begin{equation}
    \label{eq:coder}
    \left(d^Av\right)_\alpha=i\left(d\theta_\alpha-A_\alpha\right)v_\alpha.
\end{equation}
Since $d\theta_\alpha-A_\alpha=d\theta_\beta-A_\beta$ on $U_{\alpha\beta},$ this real-valued 1-form is globally defined on $\M\backslash\{b_1,\ldots,b_d\}.$ We conclude that to each such section $v$ there corresponds a one-form that can be simply written as $d\theta-A.$ 

\section{Canonical harmonic sections}
Before computing the renormalized energy, we need to introduce and describe the properties of {\em canonical harmonic sections}.

\begin{definition} Following \cite{IgnatJerrard} we define a canonical harmonic section with singularities at $b_1,\dots,b_d$ to be a smooth section of modulus $1$ on $\M\setminus\{b_1,\ldots,b_d\}$ with degree $+1$ about each of the singularities, and such that the corresponding globally defined form $\omega :=d\theta-A$ satisfies
\begin{equation}\label{canonical}d\omega = 2\pi\left(\sum_{j=1}^d\delta_{b_j}\right) d\vol - dA,\quad d^*\omega = 0.\end{equation}
\end{definition}
We now describe the set of these canonical harmonic one-forms. 
First, we have:
\begin{lem} A smooth section of modulus $1$ on $\M\setminus\{b_1,\ldots,b_d\}$ with degree $+1$ about each of the singularities exists iff $d = e(\L)$. 
\end{lem} 
\begin{proof} For the necessity of the condition $d = e(\L)$, see for instance \cite{BottTu}, Theorem 11.16. 

We sketch the proof that this condition is also sufficient. Assume $d = e(\L)$ and $b_1,\dots,b_d$ are distinct points in $\M$. Let $\psi_{b,d}$ be any solution to the equation
\begin{equation}\label{psibd}-\Delta \psi_{b,d}=2\pi\sum_{j=1}^d\delta_{b_j}-*dA,\quad \int_\M\psi_{b,d}=0.\end{equation} 
(In contrast to \cite{IgnatJerrard}, we will view $\psi_{b,d}$ as a function rather than as a two-form.)
 Since $d = e(\L)$,  the integral over $\M$ of the right-hand side is zero and therefore such a solution exists. Then $\omega = *d\psi_{b,d}$ satisfies 
 $$d\omega = 2\pi\left(\sum_{j=1}^d\delta_{b_j}\right) d\vol - dA,$$
 so that if we define $\theta_\alpha$ to be such that $d\theta_\alpha = \omega + A_\alpha$, then $e^{i\theta_\alpha}$ is well defined on $U_\alpha\setminus\{b_1,\ldots,b_d\}$ and has a singularity of degree one at each $b_k\in U_\alpha$.

Moreover, $d(\theta_\beta - \theta_\alpha - \phi_{\alpha\beta}) = 0$ for any $\alpha$, $\beta$, so that there exist constants $c_{\alpha\beta}$ such that $\theta_\beta - \theta_\alpha - \phi_{\alpha\beta} = c_{\alpha\beta}$. Note that each $\theta_\alpha$ is only defined modulo $2\pi$ due to the presence of singularities, but $\theta_\beta - \theta_\alpha$ is a well-defined function on $U_{\alpha\beta}$. 

From the definition of $c_{\alpha\beta}$, we have that 
$\varepsilon_{\alpha\beta\gamma} := c_{\alpha\beta} + c_{\beta\gamma}+c_{\gamma\alpha} $ belongs to $2\pi\mZ$ for any $\alpha$, $\beta$, $\gamma$. Therefore there exists $\tilde c_{\alpha\beta},\tilde c_{\beta\gamma},\tilde c_{\gamma\alpha}\in 2\pi\mZ$ such that 
$\varepsilon_{\alpha\beta\gamma} = {\tilde c}_{\alpha\beta} + {\tilde c}_{\beta\gamma}+{\tilde c}_{\gamma\alpha}$. 
Then $\sigma_{\alpha\beta} := c_{\alpha\beta} - {\tilde c}_{\alpha\beta}$ satisfies 
$\sigma_{\alpha\beta} + \sigma_{\beta\gamma}+\sigma_{\gamma\alpha} = 0$. Therefore there exists functions $f_\alpha:U_\alpha\to \mR$ such that
$\sigma_{\alpha\beta} = f_\beta - f_\alpha$, for any $\alpha$, $\beta$. Then we let ${\tilde\theta}_\alpha = \theta_\alpha - f_\alpha$, so that  
$${\tilde\theta}_\beta - {\tilde\theta}_\alpha - \phi_{\alpha\beta} = {\tilde c}_{\alpha\beta},$$
hence $e^{i\tilde\theta_\beta} = e^{i\tilde\theta_\alpha}e^{i\phi_{\alpha\beta}}$ and therefore $\{e^{i{\tilde\theta}_\alpha}\}$ defines a global section of the bundle.
\end{proof}

\begin{rmk}
    The section $u_0$ constructed in the above lemma has the following behavior near each point $b_k$: Choosing coordinates for which $b_k\in U_\alpha$ is the origin, we have $(u_0)_\alpha(x) =  e^{if(x)}x/|x|$ for some smooth function $f$ defined on a neighborhood of the origin.
\end{rmk}

 Next we describe the set of canonical harmonic sections. Let $u_0$ be a section of modulus 1 satisfying the degree $+1$ condition about the points $b_1,\ldots,b_d$ as given by the previous lemma and let $u$ be a canonical harmonic section. Then $u/u_0$ is a well-defined smooth function $v$ from $\M$ to the unit circle in $\mC$. 

If $\eta$ is the differential of the phase of $u/u_0$, then it is a closed one-form on $\M$. The  integral of $\eta$ along any closed path --- or, equivalently, along paths $\gamma_1,\dots,\gamma_g$ which generate $\pi_1(\M,\mZ)$ ---  belongs to  $2\pi\mZ$. We will say that $\eta$ has fluxes in $2\pi\mZ$.

It follows that the  one-form $$\omega=d\theta-A$$ on $\M\backslash\{b_1,\ldots,b_d\}$ that corresponds to $u$ is equal to $ \omega_0 + \eta$. Here $\omega_0$ is the form corresponding to $u_0$ (which is singular at the $b_k$'s) and $\eta$ is closed and has fluxes in $2\pi\mZ$. Moreover, from  \eqref{canonical} and \eqref{psibd}, we have that $\lambda := \omega - *d\psi_{b,d}$ is harmonic. Therefore, since $\eta = \lambda + *d\psi_{b,d} - \omega_0$ we deduce that  $\omega$ belongs to 
\begin{equation}\label{omegabd}\Omega_{b,d} := \left\{*d\psi_{b,d} + \lambda\mid\text{$\lambda$ is harmonic and $\lambda + *d\psi_{b,d} - \omega_0$ has integer fluxes.}\right\}.\end{equation}
Reciprocally, if $\omega=*d\psi_{b,d} + \lambda\in\Omega_{b,d}$, then $\eta = \lambda + *d\psi_{b,d} - \omega_0$ has integer fluxes and is closed from \eqref{psibd} and \eqref{canonical} using the fact that $\lambda$ is harmonic. Therefore $\eta$ is the differential of the phase of some smooth $v:\M\to\mC$ of modulus one, and we can readily check that $u = u_0 v$ is a canonical harmonic section. 

We recover the description of \cite{IgnatJerrard} by considering a set $\gamma_1,\ldots,\gamma_{2g}$ of generators of $\pi_1(\M)$ and their Poincaré duals  $\eta_1,\dots,\eta_{2g}$ which generate $H^1(\M)$, where $g$ is the genus of $\M$. If we choose $\eta_1,\dots,\eta_{2g}$ to be harmonic, then any harmonic one-form $\eta$ is of the form $\sum_k \Phi_k \eta_k + df$ for some harmonic function $f$ and from the definition of the Poincaré dual, $\int_{\gamma_\ell}\eta = \langle \eta,\eta_\ell\rangle$. Therefore we have reproved part of Theorem 2.1 in \cite{IgnatJerrard}.

\begin{prop}
\label{prop:1} Let as above $\eta_1,\dots,\eta_{2g}$ be harmonic Poincaré duals of a set of simple closed curves which generate $H^1(\M)$ and $\omega_0 = -i {u_0}^{-1} d_A u_0$ for some section $u_0$ of modulus one satisfying the degree $+1$ condition about the points $b_1,\ldots,b_d$.

Then a  one-form $\omega$ corresponds to a canonical harmonic section satisfying the degree conditions about the points $b_1,\ldots,b_d$ if and only if 
\[\omega=*d\psi_{b,d}+\sum_{k=1}^{2g}\Phi_k\eta_k+df,\]
where $f:\M\to\mathbb{R}$ is harmonic, where $\psi_{b,d}$ is defined in \eqref{psibd}, and where for each $k$ the coefficient $\Phi_k$ satisfies the condition
\begin{equation}\label{lattice}\sum_{1\le l\le 2g} a_{k\ell}\Phi_\ell-\zeta_k\in 2\pi\mathbb Z,\end{equation}
with 
\[a_{k,\ell}= \langle \eta_k,\eta_\ell\rangle,\quad \zeta_k = \langle\omega_0 - *d\psi_{b,d}, \eta_k\rangle.\]
\end{prop}

\section{Proof of Theorem 1}

Step 1 : First we compute an upper-bound for minimizers of $E_\eps$. The construction differs from that in section 9.2 of \cite{IgnatJerrard} only in that the bundle may be different than the tangent bundle of $\M$. We only sketch it. Consider points $b_1,\dots,b_d$, and a canonical harmonic section $v$ of $\L$ of modulus one with singularities of degree $1$ at each of the points. Outside $\cup_k B(b_k,\sqrt\eps)$ we let $u_\eps = v$. Near a point $b_k\in U_\alpha$, we use exponential coordinates $x$ on $\M$ centered at $b_k$. With a slight abuse of notation, we denote by $v_\alpha$ the trivialization of $v$ in the patch $U_\alpha$, expressed in exponential coordinates.

Then, following the argument in \cite{IgnatJerrard}, we deduce the  behavior for $v_\alpha$ near $b_k$:
$$v_\alpha(x) = e^{i(\theta+\eta_k(x))},$$
where $\eta_k$ is a $C^1$ real-valued function defined in a neighborhood of the origin and $\theta$ is the polar angle in $\mR^2$.

This behavior allows us to glue $v_\alpha$ defined outside $B(b_k,\sqrt\eps)$ with the map $u_0(x/\eps)$ defined inside the ball, where $u_0$ is the standard radial vortex (see for example \cite{BBH}). This will require an interpolation on the annulus $B(b_k,2\sqrt\eps)\setminus B(b_k,\sqrt\eps)$, but this layer adds only $o(1)$ to the energy of the resulting test map $v_\eps$, as $\eps\to 0$.

We have, neglecting the interpolation layer,  
$$E_\eps(v_\eps,\M\setminus\cup_k B(b_k,\sqrt\eps)) = \pi d \log\frac1{\sqrt\eps} + \tilde E[v,b_1,\dots,b_d] + o(1),$$
from \eqref{Wdefn}.
On the other hand, from the well-known properties of the radial vortex we have
$$E_\eps(v_\eps,\cup_k B(b_k,\sqrt\eps) = \pi d \log\frac{\sqrt{\eps}}{\eps} + d\gamma + o(1),$$
where $\gamma$ is defined in \cite{BBH}.
Adding these estimates, we find that 
$$E_\eps(v_\eps) = \pi d \log\frac1{\eps} + \tilde E[v,b_1,\dots,b_d] + o(1).$$
Minimizing over all canonical harmonic sections $v$  with singularities of degree one about the points $b_1,\ldots,b_d$, we find that a minimizer $u_\eps$ of $E_\eps$ must satisfy, as $\eps\to 0$, 
$$E_\eps[u_\eps]\le\pi d\log{\frac{1}{\eps}}+d\gamma+W(b_1,\ldots,b_d)+o(1),$$
for any choice of $b_1,\dots,b_d$.

\medskip
Step 2 : Compactness. 

 Next, one couples the lower bound obtained through the ball construction of \cite{Jer99} and \cite{San98}, adapted to the present manifold setting, e.g., in \cite{IgnatJerrard}, and the upper bound from Step~1. It follows that, modulo a subsequence, there exist points $a_1,\dots,a_n\in\M$ such that $u_\eps$ converges in  $H^1_\mathrm{loc}$  outside $a_1,\dots,a_n$ to some $u_*$. Moreover,
 \begin{equation}
     \label{eq:balls}
     E_\eps(u_\eps,\cup_k B(a_k,\eta)\ge \pi d\log\frac\eta\eps  - C,
 \end{equation}
for any $\eta>0$, where $C$ is independent of $\eps$ and $\eta$. It follows that $u_*$ is a section of modulus $1$ such that 
\begin{equation}\label{etabound}E_0(u_*,\M\setminus\cup_k B(a_j,\eta)):= \frac12\int_{\M\setminus\cup_k B(a_j,\eta)}|d_Au_*|^2\le \pi d\log\frac1\eta +C.\end{equation}
On the other hand, since $u_\eps$ is a critical point of $E_\eps$ we have $d^*\omega_\eps = 0$, where $\omega_\eps = iu_\eps\cdot d_Au_\eps$. Passing to the limit, we find that $\omega_* = iu_*\cdot d_Au_*$ satisfies $d^*\omega_* = 0$. But outside the points we have $|u_*| = 1$ and $\omega_* = d\theta_* - A,$ so that $d\omega_* = -dA$. Thus 
$$\text{$d^*\omega_* = 0$ and $d\omega_* = -dA$ on $\M\setminus\cup_j B(a_j,\eta)$.}$$

Step 3 : We now check what happens at the points $a_1,\dots,a_n$. From \eqref{eq:balls}, the energy inside the balls of radius $\eta$ and centered at $a_1,\dots,a_n$ is $\pi d \log\eta/\eps - C$. Then, for $u_*$ the energy in the annuli $B(a_j,\sqrt\eta)\setminus B(a_j,\eta)$ is at least $\frac12\pi\sum_j {d_j}^2|\log\eta|$ (cf. \cite{San98}). From the upper bound, taking $\eta$ small we find $\sum_j {d_j}^2 = d$, hence $n=d$ and $d_j = 1$ for every $j$. 

From the analog of \eqref{etabound} applied to the annulus $\A_\eta:=B(a_j,\sqrt\eta)\setminus B(a_j,\eta)$ we have that, for any $j$, 
\begin{equation}\label{annulusbound}E_0(u_*,B(a_j,\sqrt\eta)\setminus B(a_j,\eta))\le \pi\log\frac{\sqrt\eta}\eta + C.\end{equation}
We choose one of the points $a_j\in U_\alpha$ and let $u = (u_*)_\alpha = e^{i\theta_*}$. We choose exponential coordinates on $\M$ around $a$, and let $(r,\theta)$ be the corresponding polar coordinates, and $\omega_{a_l} = d\theta - A$. Then, using \eqref{annulusbound} and the fact that $u_*$ has degree one at $a_j$, we obtain
$$\int_{\A_\eta}|\omega_* - \omega_{a_j}|^2 = \int_{\A_\eta}|\omega_*|^2 - \int_{\A_\eta}|\omega_{a_j}|^2 +O(1)\le C.$$

Next, we choose a $C^\infty$-cutoff function $\varphi$ equal to $1$ on $B(a_k,\eta)$ and $0$ outside $B(a_k,\sqrt\eta),$ such that 
$$\int_{\A_\eta}|d\varphi|^2\le\frac C{|\log\eta|}. $$
Then, integrating by parts,
$$\int_{\M} (d\omega_* - d\omega_{a_j})\varphi = \int_{\A_\eta} (\omega_* - \omega_{a_j})\wedge d\varphi\le \|\omega_* - \omega_{a_j}\|_{L^2(\A_\eta)}\|d\varphi\|_{L^2(\A_\eta)}\le \frac C{|\log\eta|}.$$
Since in a neighborhood of $a_j$ we have $d\omega_{a_j} = 2\pi\delta_{a_k}d\vol-dA$, by letting $\eta\to 0$ we deduce that the part of $d\omega_*$ supported at the singleton $\{a_j\}$ is precisely  $2\pi\delta_{a_j}d\vol$. We conclude similarly that $\{a_j\}$ lies outside of the support of $d^*\omega_*$, and finally that 
$$\text{$d^*\omega_* = 0$ and $d\omega_* = 2\pi\left(\sum_j\delta_{a_j}\right)d\vol-dA$.}$$
Therefore $u_*$ is a canonical harmonic section. 

Step 4 : Lower bound up to $o(1)$ with the renormalized energy. Consider a subsequence, not relabeled, such that $u_\eps$ converges to $u_*$ in $C^1_{\loc}(\M\setminus\{a_1,\dots,a_d\})$. In a neighborhood of $a_j$ in a patch $U_\alpha$ we have $u_{*,\alpha} = e^{i\theta_{j,\alpha}}$, where the phase $\theta_{j,\alpha}$ is defined modulo $2\pi$ and is harmonic on $U_\alpha$. We choose exponential coordinates in a neighborhood of $a_j$, so that $a_j$ corresponds to the origin and the metric at the origin is the identity matrix.

Then, in these coordinates, $u_{*,\alpha}(x) \approx x/|x|$ on $\partial B_\rho(0)$ as $\rho\to 0$. Therefore, as in \cite{BBH} and \cite{IgnatJerrard}, choosing $\rho>0$ small enough,  
\begin{equation}\label{low1} \liminf_{\eps\to 0} \left(E_\eps(u_\eps,B_\rho(a_j)) - \pi\log\frac\rho\eps\right)\ge \gamma + o_\rho(1),\end{equation}
where $\gamma$ is defined as in \cite{BBH}. 

On the other hand, by definition of the renormalized energy in \eqref{Wdefn}, we have 
\begin{equation}\label{low2}\frac{1}{2}\int_{\M\backslash\cup_{j=1}^d B_\rho(a_j)}|\nabla_A u_*|^2\ge \tilde{E}[u_*,a_1,\ldots,a_d]+ \pi d\log\frac1{\rho} + o_\rho(1).\end{equation}
Combining \eqref{low1} and \eqref{low2}, we find that,
$$\liminf_{\eps\to 0}\left( E_\eps(u_\eps) - \pi d\log\frac1\eps\right)\ge  d\gamma + \tilde{E}[u_*,a_1,\ldots,a_d].$$
Combining with step 1, we get the desired result.

\section{The renormalized energy}
\label{sec:W}
In this section we compute a more explicit formula for the renormalized energy $W$.  
\begin{definition} We define the lattice $\LL(b_1,\dots,b_d)\subset \mR^{2g}$ to be the set of vectors $$(\Phi_1,\dots,\Phi_{2g})$$ such that \eqref{lattice} holds.
The function $G(x,y)$ is the solution to 
$$- \Delta G(\cdot,y) = \delta_y - \frac1{\vol(\M)}.$$
The regular part of $G$ is 
$$H(x,y) = G(x,y)  + \frac1{2\pi} \log(\dist_\M(x,y)),$$
where $\dist_\M$ denotes the intrinsic distance on $\M$.

The function $\psi_0$ is the solution of 
$$- \Delta \psi_0 = -\kappa + \bar \kappa,\quad \bar\kappa = \frac1{\vol(\M)}\int_\M \kappa,$$
where $\kappa = *dA$ is the curvature of the connection $A$.
\end{definition}
\begin{prop}\label{explicit}
Given any $\vec b = (b_1,b_2,\ldots,b_d)\in \M^d$, 
$$W(b_1,\ldots,b_d) = \inf_{\Phi\in\LL(\vec b)} W(b_1,\ldots,b_d,\Phi),$$
where 
\begin{align}
W(b_1,\ldots,b_d,\Phi)& = 4\pi^2\sum_{1\le i<j\le d} G(b_i,b_j)\nonumber \\ &+ 2\pi \sum_{i=1}^n \left(\pi H(b_i,b_i) +\psi_0(b_i)\right) +\frac12|\Phi|^2 +\int_\M \frac{|d\psi_0|^2}{2}.\label{Wdefn2}
\end{align}
\end{prop}
\begin{proof}
    The renormalized energy $W$ is defined in Theorem~1 as the infimum of $\tilde{E}[v,b_1,\ldots,b_d]$ given in \eqref{Wdefn}, where $v$ is a canonical harmonic section.

    These canonical harmonic sections, modulo multiplication by a constant phase, are in one-to-one correspondence with the fluxes $\Phi\in\LL(b_1,\ldots,b_d)$.

    The value of \eqref{Wdefn} for the canonical harmonic section corresponding to fluxes $\Phi$ was computed in \cite{IgnatJerrard} and is precisely  $W(b_1,\ldots,b_d,\Phi)$ defined above. Note that in this computation, only the one-form $\omega$ is involved and therefore it is irrelevant whether it is obtained from a section of a bundle $\L$ (in our case) or a tangent vector field on $\M$ (as is the case in \cite{IgnatJerrard}). The only difference is that in our case $\kappa$ is the curvature of the connection $A$ whereas it was the Gauss curvature of $\M$ in \cite{IgnatJerrard}, i.e. the curvature of the Levi-Civita connection on the tangent bundle of $\M$.  
\end{proof}

\section{Higher rank tensor bundles over surfaces}
\label{sec:high}

In this section we describe a family of complex line bundles $\L_k$, one for each Euler number $k = e(\L_k) \geq 2$, to which our previous work can be applied.  We first introduce some standard facts about tensor products of copies of $\R{2}$.  Then, we use this structure to describe the bundles $\L_k$.

Let $\{e_1, e_2\}$ be a fixed orthonormal basis for $\R{2}$ and $k \geq 2$ be an integer. We introduce $k$-tuples of vectors
\[e_{i_1},\ldots,e_{i_k}\]
drawn from the set $\{e_1, e_2\}$, where $i_j\in\{1,2\}$ for every $j\in\{1,\ldots,k\}.$ Considering all choices $\mu=\left\{i_1,\ldots,i_k\right\}$ of distinct arrangements of 1's and 2's, it is easy to show that the collection
\begin{equation}
    \label{eq:fmu}
    f_\mu:=e_{i_1}\otimes\ldots\otimes e_{i_k}
\end{equation}
forms a basis for $\otimes^k \R{2}:=\underbrace{\R{2}\otimes ... \otimes \R{2}}_\text{k \,\,\,\mbox{times}}.$ Further, the inner product of two elements of $\otimes^k \R{2}$ can be obtained by extending the expression 
\begin{equation}
\left\langle f_{\mu_1},f_{\mu_2}\right\rangle:=\Pi_{j=1}^k\left\langle e_{i_j^1},e_{i_j^2}\right\rangle \label{inner_prod}
\end{equation}
for the inner product between the elements of the basis.

For $\sigma \in S_k$, where $S_k$ is the group of permutations in $k$ indices, define the permutation operator
$$
T_\sigma : \otimes^k \R{2} \to \otimes^k \R{2}
$$
by extending by linearity to the entire $\otimes^k \R{2}$ the action
$$
T_\sigma(f_\mu) = e_{i_{\sigma(1)}}\otimes ... \otimes e_{i_{\sigma(k)}}
$$
of $T_\sigma$ on the elements of the basis of $\otimes^k \R{2}$. Then we can define
\begin{equation}\label{def:sym_tensors}
\otimes^k_{sym} \R{2} = \{\mathcal{Q}\in \otimes^k\R{2} : \,\,\,  T_\sigma(\mathcal{Q}) = \mathcal{Q} \,\,\, \mbox{for all}\,\,\,  \sigma \in S_k\},
\end{equation}
which will be referred to as the set of symmetric $k$-tensors on $\R{2}$.
A basis for $\otimes^k_{sym}\R{2}$ is the collection
\begin{equation}\label{def:basis_sym_tensors}
\sum_{\sigma} f_{\sigma(\mu)},
\end{equation}
for all {\it non-decreasing} $k$-tuples $\mu$ with values in $\{1,2\}$. One can identify each element of this basis through the number of 1's in $\mu,$ hence for each $j\in\{0,\ldots,k\},$ we can define
$$h(k,j)=\sum_{\sigma} f_{\sigma(\mu)},$$
where $\mu$ is the non-decreasing $k$-tuple with values in $\{1,2\}$ that begins with exactly $j$ entries equal to $1$. It immediately follows that $\otimes^k_{sym}\R{2}$ has dimension $k+1$.


We now define the trace operator
$$
{\rm{tr}}_{12} : \otimes^k\R{2} \to \otimes^{k-2}\R{2},
$$
which is the contraction with respect to the first two indices: 
$$
{\rm{tr}}_{12}(v_1\otimes ... \otimes v_k) = \left\{\begin{array}{ll}\langle v_1, v_2\rangle \left(v_3\otimes ... \otimes v_k\right), & k> 2,\\\langle v_1, v_2\rangle, & k=2.\end{array}\right.
$$
Note that for $k \geq 2$, if $\mathcal{Q}\in \otimes^k_{sym}\R{2}$, then ${\rm{tr_{12}}}(\mathcal{Q}) \in \otimes^{k-2}_{sym}\R{2}$.  We will also define the space of symmetric, traceless tensors in $\R{2}$ as follows:
\begin{equation}\label{def:sym_tr_tensors}
\otimes^k_{sym, tr} \R{2} = \{ \mathcal{Q}\in \otimes^k_{sym} : {\rm tr}_{12}(\mathcal{Q}) = 0\}.
\end{equation}

\medskip
\medskip

Furthermore, relative to the bases $\{h(k, j)\}_{j=0}^k$ and $\{h(k-2, j)\}_{j=0}^{k-2}$, for $\otimes^k_{sym}\R{2}$ and $\otimes^{k-2}_{sym}\R{2},$ respectively, the $(k-1)\times (k+1)$ matrix corresponding to the mapping $\rm{tr}_{12}$ is ${\left(\rm{tr}_{12}\right)}_{ij}=\delta_{i,j}+\delta_{i,j-2}.$ Note that the operator $\rm{tr}_{12} : \otimes^k_{sym}\R{2} \to \otimes^{k-2}_{sym}\R{2}$ is onto and so by the Kernel-Image theorem, $\rm{ker}(\rm{tr}_{12})$ has dimension  $2.$ Hence, it is easy to verify that a basis for this space is given by two tensors
\begin{equation}
\mathcal{Q}_o(k) = \sum_{1 \leq 2i+1 \leq k} (-1)^i h(k, 2i+1), \,\,\,  \mathcal{Q}_e(k) = \sum_{0 \leq 2i \leq k} (-1)^i h(k, 2i).\label{basis_sym_trace}
\end{equation}
Further, when $k \geq 2$, one can check that
\begin{align}
\mathcal{Q}_o(k+1) &= e_1 \otimes \mathcal{Q}_e(k) + e_2 \otimes \mathcal{Q}_o(k), \nonumber\\
\mathcal{Q}_e(k+1) &= -e_1 \otimes \mathcal{Q}_o(k) + e_2 \otimes \mathcal{Q}_e(k), \label{tensor_growth}
\end{align}
and
$$
\langle \mathcal{Q}_o , \mathcal{Q}_e\rangle = 0, \,\,\, \langle \mathcal{Q}_o , \mathcal{Q}_o\rangle = \langle \mathcal{Q}_e , \mathcal{Q}_e\rangle.
$$
In terms of our discussion so far, we define the set $\otimes^k_{sym,tr} \R{2}$ of rank-$k$, symmetric, traceless tensors in $\R{2}$ as the kernel of the operator $\rm{tr}_{12}$ in $\otimes^k_{sym} \R{2}$.

Next, given $0 \leq \alpha < 2\pi$, define the rotation matrix 
$$
R_\alpha = \left (  \begin{array}{cc} \cos{\alpha} & -\sin{\alpha} \\ \sin{\alpha} & \cos{\alpha} \end{array}\right )
$$
by $\alpha$ in $\R{2}$. Recall that for a linear map $T: \R{2} \to \R{2}$, we can define $\otimes^k T : \otimes^k \R{2} \to \otimes^k \R{2}$ by the condition
$$
(\otimes^k T)(v_1\otimes ... \otimes v_k) = (Tv_1) \otimes ... \otimes (Tv_k).
$$
We have the following lemma.
\begin{lem} The action of $\otimes^k R_\alpha$ on $\rm{ker}(\rm{tr}_{12})$ is determined by
$$
(\otimes^k R_\alpha)(\mathcal{Q}_o(k)) = \cos(k\alpha) \mathcal{Q}_o(k) + \sin(k\alpha) \mathcal{Q}_e(k)
$$
and
$$
(\otimes^k R_\alpha)(\mathcal{Q}_e(k)) = -\sin(k\alpha) \mathcal{Q}_o(k) + \cos(k\alpha) \mathcal{Q}_e(k).
$$
\end{lem}
\begin{proof}
The proof is straightforward for the case $k=2$.  For $k \geq 2$, the conclusion follows by induction using equations \eqref{tensor_growth}.

\end{proof}
We remark that when $k=2,$ the set $\otimes^2_{sym, tr}\R{2}$ is precisely the set of $Q$-tensors in $\R{2}$. 

\medskip

Next, let $\M$ be a closed, compact, orientable surface embedded in $\R{3}$.  We denote by $\nu$ the globally defined normal to $\M$.  We will use the usual notation $T\M$ and $T_p\M$ for the tangent bundle of $\M$ and the tangent to $\M$ at $p \in \M$, respectively.  It is a standard fact that the inner product and Levi-Civita connection from $\R{3}$ induce a Riemannian metric, and corresponding Levi-Civita connection, in $T\M$.  This turns $T\M$ into a Riemannian vector bundle with a metric connection, which we denote by $\nabla$.

We assume throughout this section the existence of a simple cover $\{U_\alpha\}$ of $\M$ such that for each $\alpha$ there is a pair of orthonormal sections $e_j^\alpha : U_\alpha \to T\M$, $j=1, 2$ and that $\{e_1^\alpha, e_2^\alpha, \nu\}$ is a positively oriented basis in $\R{3}$.
Since $\nabla$ is the Levi-Civita connection for $\M,$ there exists a $1$-form $A_\alpha$, over $U_\alpha$ that satisfies
$$
\nabla_X e_1^\alpha = A_\alpha(X) e_2^\alpha, \,\,\, \nabla_X e_2^\alpha = -A_\alpha(X)e_1^\alpha.
$$
Recall that while $A_\alpha$ is not globally defined, $dA_\alpha$ is a globally defined $2$-form.

For $k \geq 2$ we denote by $T^k\M$ the vector bundle over $\M$ with a fiber $\otimes^k T_p\M$ over $p \in\M$.  Note that the bundle projection $\Pi_k : T^k \M \to \M$ can be defined in the obvious way for every $k \geq 2$.  It is standard to extend the Riemannian metric and connection from $T\M$ to $T^k\M$, so that $T^k\M$ is also a Riemannian vector bundle with a metric connection.  Furthermore, it is routine to check that the definitions of $\otimes^k_{sym} \R{2}$ and $\otimes^k_{sym,tr} \R{2}$ extend naturally to define the bundles $T^k_{sym}\M$ of symmetric rank-$k$ tensors and $T^k_{sym,tr} \M$ of symmetric, traceless, rank-$k$ tensors over $\M$, respectively.  These are also Riemannian bundles with a metric connection, and we will take $\L_k = T^k_{sym,tr} \M$.  These constructions are standard, and the remaining work is devoted to showing that $\L_k$ is a complex line bundle, or equivalently, a rank-2 real vector bundle, and to computing the Euler number of $\L_k$.

If $\{U_\alpha\}$ is a simple cover of $\M$, then for each $\alpha$ there is a pair of orthonormal sections $e_j^\alpha : U_\alpha \to T\M$, $j=1, 2.$ Generalizing \eqref{eq:fmu}, we introduce maps
$f^\alpha_\mu: U_\alpha \to T^k\M,$ where $f^\alpha_\mu$ denote a set of orthonormal sections over $U_\alpha$ with values in $T^k\M$. At each $p \in \M$, these sections span the fiber of $T^k \M$ over $p$. This allows us to define local trivializations
$$
\Phi_\alpha : \Pi_k^{-1}(U_\alpha) \to U_\alpha \times \otimes^k \R{2}
$$
by setting
$$
\Phi_\alpha(p, f^\alpha_\mu) = (p, f_\mu).
$$
If $U_\alpha \cap U_\beta \neq \emptyset$, as before, there is a transition function $\varphi_{\alpha\beta} \to \R{}$ such that
$$
e_1^\beta=\cos{\varphi_{\alpha\beta}} e_1^\alpha + \sin{\varphi_{\alpha\beta}} e_2^\alpha \,\,\, \mbox{and} \,\,\, e_2^\beta=-\sin{\varphi_{\alpha\beta}} e_1^\alpha + \cos{\varphi_{\alpha\beta}} e_2^\alpha.
$$
Let
$$
R_{\alpha\beta} = R_{\varphi_{\alpha\beta}}= \left (  \begin{array}{cc} \cos{\varphi_{\alpha\beta}} & -\sin{\varphi_{\alpha\beta}} \\ \sin{\varphi_{\alpha\beta}} & \cos{\varphi_{\alpha\beta}}
\end{array}\right ).
$$
For the trivializations above, the transition maps
$$
\Phi_\beta \circ \Phi_\alpha^{-1} : (U_\alpha \cap U_\beta) \times \otimes^k \R{2} \to (U_\alpha \cap U_\beta) \times \otimes^k \R{2}
$$
are given by
$$
(\Phi_\beta \circ \Phi_\alpha^{-1})(p, \mathcal{Q}) = (p, (\otimes^k R_{\alpha\beta})(\mathcal{Q})).
$$
Next, we observe $T^k\M$ becomes a Riemannian vector bundle if we endow its fibers with the inner product defined in \eqref{inner_prod}.  Further, it is standard to use Leibnitz rule to furnish $T^k\M$ with a metric connection, that we still denote by $\nabla$.  This defines $T^k\M$.

\medskip
\medskip

We now point out that the the constructions that yielded the spaces $\otimes^k_{sym} \R{2}$ and $\otimes^k_{sym, tr} \R{2}$, defined in \eqref{def:sym_tensors} and \eqref{def:sym_tr_tensors} respectively, can be applied to $T^k\M$ to define the bundles $T^k_{sym}\M$ and $T^k_{sym, tr}\M$.  It is easy to check that we can use the orthonormal sections $e_j^\alpha : U_\alpha \to T\M$, $j=1, 2$ to define constant multiples of orthonormal sections
$$
\sum_{\sigma} f_{\sigma(\mu)}^\alpha
$$
and
$$
\mathcal{Q}_o^\alpha, \mathcal{Q}_e^\alpha
$$
in $U_\alpha$, with values in $T^k_{sym}\M$ and $T^k_{sym, tr}\M$ respectively.  It is straightforward to define trivializations for $T^k_{sym}\M$ and $T^k_{sym, tr}\M$ over $U_\alpha$, starting from the sections we just described.  From the constructions in \eqref{def:sym_tensors} and \eqref{def:sym_tr_tensors}, we conclude that $T^k_{sym}\M$ is a rank-$(k+1)$ real vector bundle over $\M$, whereas $T^k_{sym, tr}\M$ is a rank-$2$ vector bundle over $\M$, or equivalently, a complex line bundle.

It is a matter of a simple calculation to show that the connection $\nabla$ from $T^k\M$ can be restricted to $T^k_{sym}\M$ to yield a metric connection in $T^k_{sym}\M$, and the same is true for $T^k_{sym, tr}\M$. We can now state the main result of this section.

\begin{lem}
\label{lem:main}
For each $U_\alpha$, the tensors $\mathcal{Q}_o^\alpha(k)$ and $\mathcal{Q}_e^\alpha(k)$ given by \eqref{basis_sym_trace} form an orthogonal frame for $T^k_{sym, tr}\M$ in $U_\alpha$, and they have the same constant length, that depends only on $k$, and not on $\alpha$. Furthermore, in each $U_\alpha,$ we also have
$$
\nabla_X \mathcal{Q}_o^\alpha = k \,A_\alpha(X)\,\mathcal{Q}_e^\alpha, \,\,\, \nabla_X \mathcal{Q}_e^\alpha = -k \,A_\alpha(X)\,\mathcal{Q}_o^\alpha.
$$
\end{lem}
\begin{proof}
Most of the Lemma is already proved in the previous section.  The claim about the covariant derivatives of $\mathcal{Q}_o$, $\mathcal{Q}_e$ follows from \eqref{tensor_growth} by induction in $k$.
\end{proof}
The final proposition is a corollary of Lemma \eqref{lem:main}.
\begin{prop}
\label{prop:porp}
    The Euler number of the bundle of symmetric, traceless, rank-$k$ tensors over $\M$ is given by
    $$ e(T^k_{sym,tr}\M) = 2k e(T\M).$$
\end{prop}
\begin{proof}
From Lemma~\eqref{lem:main}, and since  $\mathcal{Q}_o$, $\mathcal{Q}_e$ is (a constant multiple of) an orthonormal frame, the curvature of $\nabla_X$ on $T^k_{sym,tr}\M$ is $k$ times that of $\nabla_X$ on $T\M$. We conclude by integrating curvature over $\M$ (cf. \cite{MR3585539}).
\end{proof}
When $\M=\mathbb S^2$ we immediately obtain
\begin{cor}
    The Euler number of $T^k_{sym,tr}\mathbb S^2$ is  $2k$.
\end{cor}

Before concluding this section, we remark that alternatively, if one views the tangent planes $T\M_p$ at different points $p\in\M$ as one-dimensional complex lines rather than two-dimensional real planes, then this provides a different but equivalent description of the line bundles $T^k_{sym,tr}\M$, \cite{Hall}. However, because our work is inspired by applications in soft matter and meshing, we choose here to present our work in the context of real vector spaces rather than complex. 

\section{Proof of Theorem 2}

Here we consider the case where $\M=\mathbb{S}^2$, where $\L$ is the rank-$k,$ symmetric, traceless tensor bundle over $\mathbb S^2$, and where the connection $A$ is the Levi-Civita connection on $\L$. The curvature of $A$ is then constant $\kappa = 2k$.

By Proposition \eqref{prop:porp}, the Euler number of $\L$ is equal to $2k$ and, since $\mathbb S^2$ is simply connected, there is a unique canonical harmonic section associated to singularities $(b_1,\dots,b_{2k})$. Furthermore, the lattice $\LL(b_1,\dots,b_{2k})$, being of dimension $2g$, degenerates to $\{0\}$ since the genus of the sphere is $g=0$. 

Moreover,  the Green's function is given by
$$G(x,y) =  -\frac1{2\pi}\log(\abs{x-y}),$$
where $\abs{x-y}$ is the Euclidean distance between two points in $\mathbb{R}^3$. As a consequence of these observations, formula \eqref{Wdefn} reduces in this case to
$$W(a) = -2\pi\sum_{1\le i<j\le 4}\log|a_i-a_j|.$$

Then Theorem \ref{tetra} follows from \cite{MR2257398}, Theorem 1.2, since the logarithm is a completely monotonic function.
\section*{Acknowledgements}
 The research of D.G. was supported by an NSF grant DMS 2106551. E.S. thanks  Indiana University Bloomington, where part of this research was conducted, for its hospitality. The research of P.S. was supported by a Simons Collaboration grant 585520 and an NSF grant DMS 2106516.
\bibliographystyle{abbrv}

\bibliography{vescicles}
\end{document}